\documentclass[11pt,reqno]{amsart}
\usepackage{amssymb,amscd,amsbsy}
\usepackage{amssymb,amscd,amsbsy,mathrsfs}
\newcommand{\g}{\gamma}

\newcommand{\e}{\varepsilon}
\newcommand{\vk}{\varkappa}

\renewcommand{\l}{\lambda}

\newcommand{\s}{\sigma}
\renewcommand{\t}{\tau}

\newcommand{\f}{\varphi}

\newcommand{\E}{{\mathscr E}}

\newcommand{\F}{{\mathscr F}}

\newcommand{\h}{{\mathscr H}}

\newcommand{\X}{{\mathscr X}}
\newcommand{\Y}{{\mathscr Y}}

\newcommand{\C}{{\Bbb C}}

\newcommand{\R}{{\Bbb R}}
\newcommand{\Z}{{\Bbb Z}}

\newcommand{\0}{{\boldsymbol{0}}}

\newcommand{\bs}{\boldsymbol}

\newcommand{\bS}{{\boldsymbol S}}

\newcommand{\rf}[1]{(\ref{#1})}

\newcommand{\df}{\stackrel{\mathrm{def}}{=}}

\newcommand{\supp}{\operatorname{supp}}

\newcommand{\const}{\operatorname{const}}

\newcommand{\eeq}{\end{equation}}
\newcommand{\beq}{\begin{equation}}
\newcommand{\bay}{\begin{eqnarray}}
\newcommand{\ba}{\begin{align*}}
\newcommand{\ea}{\end{align*}}
\newcommand{\ey}{\end{eqnarray}}
\newcommand{\bey}{\begin{eqnarray*}}
\newcommand{\eey}{\end{eqnarray*}}

\newcommand{\be}{\infty}

\newcommand{\bl}{\blacksquare}

\newcommand{\Pf}{{\bf Proof. }}
\newcommand{\im}{\operatorname{Im}}

\newcommand{\ov}{\overline}

\newtheorem{thm}{\hspace{\parindent}Theorem}[section]

\newtheorem{lem}[thm]{\hspace{\parindent}Lemma}

\theoremstyle{remark}

\newtheorem*{rem*}{Remark}

\newcommand\Li{{\rm Lip}}

\begin{document}

\newcommand{\vse}{\vspace{.2in}}
\numberwithin{equation}{section}

\title{Functions of triples of noncommuting self-adjoint operators under perturbations of class $\bs{S_p}$}
\author{V.V. Peller}
\thanks{the author is partially supported by NSF grant DMS 1300924
and by the Ministry of Education and Science of the Russian Federation (the Agreement number N$^{\underline\circ}$ 02.à03.21.0008).}

\begin{abstract}
In this paper we study properties of functions of triples of not necessarily commuting self-adjoint operators. The main result of the paper shows that unlike in the case of functions of pairs of self-adjoint operators there is no Lipschitz type estimates in any Schatten--von Neumann norm $\bS_p$, $1\le p\le\be$, for arbitrary functions in the Besov class $B_{\be,1}^1(\R^3)$. In other words, we prove that 
for $p\in[1,\be]$,
there is no constant $K>0$ such that the inequality
\begin{align*}
\|f(A_1,B_1,C_1)&-f(A_2,B_2,C_2)\|_{\bS_p}\\[.1cm]
&\le K\|f\|_{B_{\be,1}^1}
\max\big\{\|A_1-A_2\|_{\bS_p},\|B_1-B_2\|_{\bS_p},\|C_1-C_2\|_{\bS_p}\big\}
\end{align*}
holds for an arbitrary function $f$ in $B_{\be,1}^1(\R^3)$ and for arbitrary finite rank self-adjoint operators $A_1,\,B_1,\,C_1,\,A_2,\,B_2$ and $C_2$.
\end{abstract}

\maketitle


\

%
%
%
%

\setcounter{section}{0}
\section{\bf Introduction}
\setcounter{equation}{0}
\label{In}

\

The spectral theorem for commuting self-adjoint operators implies that for commuting 
self-adjoint operators $A_1$ and $A_2$ and for a Lipschitz function $f$ on the real line
$\R$ the following Lipschitz type estimate holds
$$
\|f(A_1)-f(A_2)\|\le\|f'\|_{L^\be(\R)}\|A_1-A_2\|.
$$
The same inequality holds for the norms in Schatten--von Neumann classes $\bS_p$ with
$p\ge1$. However, for noncommuting self-adjoint operators, the situation is quite different.
A Lipschitz function $f$ on $\R$ does not have to be {\it operator Lipschitz}, i.e., the inequality
$$
|f(x_1)-f(x_2)|\le\const|x_1-x_2|,\quad x_1,~x_2\in\R,
$$
does not imply that
$$
\|f(A_1)-f(A_2)\|\le\const\|A_1-A_2\|
$$
for self-adjoint operators $A_1$ and $A_2$.
This was shown by Farforovskaya in \cite{F1}. She also proved in \cite{F2} that
there exist a Lipschitz function $f$ on $\R$ and self-adjoint operators $A_1$ and $A_2$ such that $A_1-A_2$ belongs to trace class $\bS_1$, but $f(A_1)-f(A_2)\not\in\bS_1$.

Recall that a function $f$ on $\R$ is operator Lipschitz if and only if it takes trace class perturbations to trace class increments, i.e.,
$$
A=A^*,\quad B=B^*,\quad A-B\in\bS_1\quad\Longrightarrow\quad f(A)-f(A)\in\bS_1
$$
if we admit not necessarily bounded self-adjoint operators $A$ and $B$, see \cite{AP}.

It was shown later in \cite{Mc} and \cite{Ka} that the function $x\mapsto|x|$ is not operator Lipschitz. Necessary conditions for operator Lipschitzness were obtained in \cite{Pe2} and \cite{Pe3}. In particular, it was proved in \cite{Pe2} that operator Lipschitz functions on $\R$ must belong locally to the Besov class $B_{1,1}^1(\R)$.
Note that in \cite{Pe3} stronger necessary conditions were also found.
Those necessary conditions were deduced from the trace class criterion for Hankel operators, see \cite{Pe1} and \cite{Pe4}.

On the other hand, it was proved in \cite{Pe2} and \cite{Pe3} that functions in the Besov class $B_{\be,1}^1(\R)$ are necessarily operator Lipschitz.
This result was generalized in \cite{APPS} to functions of normal operators. It was shown in \cite{APPS} that if $f$ is a function of two variables that belongs to the Besov class  
$B_{\be,1}^1(\R^2)$, then $f$ is an {\it operator Lipschitz function on $\R^2$}, i.e.,
$$
\|f(N_1)-f(N_2)\|\le\const\|f\|_{B_{\be,1}^1}\|N_1-N_2\|
$$
for arbitrary normal operators $N_1$ and $N_2$. The same Lipschitz type inequality holds in the Schatten--von Neumann norm $\bS_p$ for $p\ge1$.

This result was generalized in \cite{NP} to the case of functions of $d$-tuples of commuting self-adjoint operators: if $f$ belongs to the Besov class $B_{\be,1}^1(\R^d)$ and $(A_1,\cdots,A_d)$ and $(B_1,\cdots,B_d)$ are $d$-tuples of commuting self-adjoint operators, then
$$
\|f(A_1,\cdots,A_d)-f(B_1,\cdots,B_d)\|\le\const\|f\|_{B_{\be,1}^1}
\max_{1\le j\le d}\|A_j-B_j\|
$$
and the same inequality also holds for Schatten--von Nemann norms $\bS_p$ with $p\ge1$.

Let me also mention that in \cite{KPSS} it was shown that for an arbitrary Lipschitz function $f$ on $\R^d$ and for $p\in(1,\be)$ the following Lipschitz type inequality holds:
$$
\|f(A_1,\cdots,A_d)-f(B_1,\cdots,B_d)\|_{\bS_p}\le\const\|f\|_{\Li}
\max_{1\le j\le d}\|A_j-B_j\|_{\bS_p}
$$
for arbitrary $d$-tuples of commuting self-adjoint operators $(A_1,\cdots,A_d)$ and $(B_1,\cdots,B_d)$; the constant on the right-hand side depends on $p$. Note that 
earlier in the case $d=1$ this was established in \cite{PS}.

We refer the reader to the recent survey article \cite{AP}, which is a comprehensive study of operator Lipschitz functions. 

The  behavior of functions of pairs of {\it noncommuting} self-adjoint operators under perturbation was studied in \cite{ANP}. For a pair $(A,B)$ of not necessarily commuting self-adjoint operators the functions $f(A,B)$ are defined as double operator integrals:
$$
f(A,B)=\iint f(x,y)\,dE_A(x)\,dE_B(y)
$$
under the assumption that the double operator integral makes sense. Here  $E_A$ and $E_B$ stand for the spectral measures of $A$ and $B$.

In the case when $A$ and $B$ are finite rank self-adjoint operators (or, more general, if $A$ and $B$ have finite spectra), the operator $f(A,B)$ is defined for all functions $f$ on $\R^2$:
$$
f(A,B)=\sum_{j,k}f(\l_j,\mu_k)P_jQ_k,
$$
where 
$$
A=\sum_j\l_jP_j\quad\mbox{and}\quad B=\sum_k\mu_kQ_k
$$
are the spectral expansions of $A$ and $B$.

It turned out that the situation in the case of noncommuting operators is different. It was shown in \cite{ANP} that if $f$ belongs to the Besov class $B_{\be,1}^1(\R^2)$ and $1\le p\le2$, then, as in the case of commuting operators, the following Lipschitz type estimate holds:
$$
\|f(A_1,B_1)-f(A_2,B_2)\|_{\bS_p}\le\const\|f\|_{B_{\be,1}^1}
\max\big\{\|A_1-A_2\|_{\bS_p},\|B_1-B_2\|_{\bS_p}\big\}
$$
for arbitrary pairs $(A_1,B_1)$ and $(A_2,B_2)$ of not necessarily commuting self-adjoint operators.

On the other hand, it was shown in \cite{ANP} that unlike in the case of commuting operators there is no Lipschitz type estimate in the norm of $\bS_p$ for $p>2$ as well as in the operator norm. In other words, if $p>2$, there is no constant $K$ such that
$$
\big\|f(A_1,B_1)-f(A_2,B_2)\big\|_{\bS_p}\le K\|f\|_{B_{\be,1}^1}
\max\big\{\|A_1-A_2\|_{\bS_p},\|B_1-B_2\|_{\bS_p}\big\}
$$
for arbitrary finite rank self-adjoint operators $A_1$,$B_1$, $A_2$ and $B_2$.
The same is true in the operator norm.

In this paper we deal with functions of triples of not necessarily commuting self-adjoint operators. For a triple $(A,B,C)$ of not necessarily commuting self-adjoint operators and a function $f$ on $\R^3$, the operator $f(A,B,C)$ is defined as the triple operator integral
$$
f(A,B,C)=\iiint f(x,y,z)\,dE_A(x)\,dE_B(y)\,dE_C(z)
$$
in the case when the triple operator integral is defined. Again, if $A$, $B$ and $C$ have finite spectra, the triple operator integral on the right is well defined for all functions $f$ and
$$
f(A,B,C)=\sum_{\l\in\s(A),\,\mu\in\s(B),\,\nu\in\s(C)}f(\l,\mu,\nu)
E_A(\{\l\})E_B(\{\mu\})E_C(\{\nu\}).
$$

The main objective of this paper is to show that {\it unlike in the case of functions of two noncommuting self-adjoint operators}, there is no Lipschitz type estimate in the norm od $\bS_p$, $1\le p\le\be$,
for functions in the Besov class $B_{\be,1}^1(\R^3)$. In other words, there is no constant $K>0$ such that
\begin{align*}
\big\|f(A_1,B_1,C_1)&-f(A_2,B_2,C_2)\big\|_{\bS_1}\\[.2cm]
&\le K\|f\|_{B_{\be,1}^1}
\max\{\|A_1-A_2\|_{\bS_1},\|B_1-B_2\|_{\bS_1},\|C_1-C_2\|_{\bS_1}\}
\end{align*}
for arbitrary functions $f$ in $B_{\be,1}^1(\R^3)$ and arbitrary finite rank self-adjoint operators $A_1$, $B_1$, $C_1$, $A_2$, $B_2$ and $C_2$.
In the special case $p=1$ a different poof was given in \cite{Pe7}. Note, however, that the method used in \cite{Pe7} does not work in the case $p=2$.

The main result of this paper terminates the chain of the results of the papers \cite{Pe2} and \cite{Pe3} (with Lipschitz type estimates in the operator norm and the trace norm for self-adjoint operators and functions of Besov class $B_{\be,1}^1(\R)$), \cite{APPS} (Lipschitz type estimates in the operator norm and the trace norm for normal operators
and functions of class $B_{\be,1}^1(\R^2)$), \cite{NP}
(Lipschitz type estimates in the operator norm and the trace norm
for $d$-tuples ofcommuting 
self-adjoint operators and functions of class 
$B_{\be,1}^1(\R^d)$) and, finally, \cite{ANP} (Lipschitz type estimates in the Schatten--von Neumann norms $\bS_p$, $1\le p\le2$, for  pairs of noncommuting self-adjoint operators and functions of class
$B_{\be,1}^1(\R^2)$). The results of \S\:\ref{osnre} of this paper
show that as soon as we admit three noncommuting self-adjoint operators, it becomes impossible to obtain such Lipschitz type estimates for arbitrary functions of class 
$B_{\be,1}^1(\R^3)$ in the norm of $\bS_p$ for any $p\in[1,\be]$.

� \S\:\ref{dois} of this paper we collect necessary information on multiple operator integrals, while in \S\:\ref{besov} we define theBesov classes 
$B_{\be,1}^1(\R^d)$ and briefly describe their properties.

\

\section{\bf Multiple operator integrals}
\setcounter{equation}{0}
\label{dois}

\

Double operator integrals appeared in the paper \cite{DK} by Daletskii and S.G. Krein. 
Later the beautiful theory of double operator integrals was created by Birman and Solomyak in \cite{BS1}, \cite{BS2} and \cite{BS3}. 

Let $(\X,E_1)$ and $(\Y,E_2)$ be spaces with spectral measures $E_1$ and $E_2$
on a Hilbert space $\h$, let $T$ be a bounded linear operator on $\h$ and let $\Phi$ be a bounded measurable function on $\X\times\Y$. {\it Double operator integrals} are expressions of the form
\bay
\label{dvoi}
\int\limits_\X\int\limits_\Y\Phi(x,y)\,d E_1(x)T\,dE_2(y).
\ey
Birman and Solomyak's starting point is the case when $T$ belongs to the Hilbert--Schmidt class $\bS_2$. In this case they defined double operator integrals of the form \rf{dvoi}
for arbitrary bounded measurable $\Phi$ and proved that
$$
\left\|\int\limits_\X\int\limits_\Y\Phi(x,y)\,d E_1(x)T\,dE_2(y)\right\|_{\bS_2}
\le\|\Phi\|_{L^\be}\|T\|_{\bS_2}
$$
(see \cite{BS1}).

To define double operator integrals for arbitrary bounded linear operators $T$ in the general case, restrictions on $\Phi$ must be imposed. Double operator integrals for arbitrary bounded operators $T$ can be defined for functions $\Phi$ that are {\it Schur multipliers} with respect to the spectral measures $E_1$ and $E_2$, see \cite{BS1}, \cite{Pe2}, \cite{Pi} and \cite{AP} for details.

In this paper we need double operator integrals only in the case when the spectral measures $E_1$ and $E_2$ are atomic and have finitely many atoms. 
We say that a spectral measure $E$ on a set $\X$ {\it is atomic and has finitely many atoms} if all subsets of $\X$ are measurable and there are points $a_1,\cdots,a_n$ in $\X$, called the {\it atoms}, such that 
$$
E\left(\X\setminus\bigcup_{j=1}^na_j\right)=\0\quad\mbox{and}\quad E(\{a_j\})\ne\0,\quad1\le j\le n.
$$

In the case when the spectral measures $E_1$ and $E_2$ are atomic with finitely many atoms, we
can define double operator integrals of the form \rf{dvoi} for arbitrary functions $\Phi$ by
\bay
\label{konato}
\int\limits_\X\int\limits_\Y\Phi(x,y)\,d E_1(x)T\,dE_2(y)=
\sum_{j,k}\Phi(a_j,b_k)E_1(\{a_j\})TE_2(\{b_k\}),
\ey
where the $a_j$ and the $b_k$ are the atoms of $E_1$ and $E_2$.

Under these assumptions, the norm of the linear transformer
$$
T\mapsto\iint\Phi(x,y)\,d E_1(x)T\,dE_2(y)
$$
(both in the operator norm and in the trace norm) is equal to the norm of the matrix
$\{\Phi(a_j,b_k)\}$ in the space of matrix Schur multipliers, i.e., the norm of the matrix transformer
$$
\{\g_{jk}\}\mapsto\{\Phi(a_j,b_k)\g_{jk}\}
$$
in the operator norm (or in the trace norm), see \cite{AP}.

Double operator integrals play an important role in perturbation theory. In particular, a special role is played by the following formula:
\bay
\label{razdrai}
f(A)-f(B)=\iint\limits_{\R\times\R}\frac{f(x)-f(y)}{x-y}\,dE_A(x)(A-B)\,dE_B(y),
\ey
which holds for arbitrary self-adjoint operators $A$ and $B$ with bounded
$A-B$ and for arbitrary operator Lipschitz functions $f$ on $\R$, see \cite{BS3}
and \cite{AP}.

In this paper we consider only operators with finite spectra, in which case formula 
\rf{razdrai} holds for arbitrary functions $f$ on $\R$; moreover, the divided difference
$$
(x,y)\mapsto\frac{f(x)-f(y)}{(x-y)}
$$
can be extended to the diagonal $\{(x,x):~x\in\R\}$
arbitrarily, i.e., the values of the divided difference on the diagonal do not affect the right-hand side of \rf{razdrai}. This can be verified elementarily.

{\it Multiple operator integrals}
$$
\underbrace{\int\cdots\int}_m\Phi(x_1,\cdots,x_m)
\,dE_1(x_1)T_1\,dE_2(x_2)T_2\cdots\,dE_{m-1}(x_{m-1})T_{m-1}\,dE_m(x_m)
$$
were defined for functions $\Phi$ in the (integral) {\it projective tensor product}
of the spaces $L^\be(E_j)$,
$j=1,\cdots,n$, in \cite{Pe5}. Later multiple operator integrals were defined in \cite{JTT} for functions $\Phi$ in the Haagerup tensor products of $L^\be$ spaces. We refer the reader to the survey article \cite{Pe6} for detailed information about multiple operator integrals.

Again, in this paper we consider only atomic spectral measures with finitely many atoms, in which case multiple operator integrals can be defined for arbitrary functions $\Phi$. Indeed, consider for simplicity the case of triple operator integrals. Suppose that $a_j$, $B_k$ and $c_l$ are the atoms of $E_1$, $E_2$ and $E_3$ and $\Phi$ is an arbitrary function. Then
\begin{align*}
\iiint\Phi(x_1,x_2,x_3)&\,dE_1(x_1)T_1\,dE_2(x_2)T_2\,dE_3(x_3)\\[.2cm]
&\df
\sum_{j,k,l}\Phi(a_j,b_k,c_l)E_1(\{a_j\})T_1E_2(\{b_k\})T_2E_3(\{c_l)\}.
\end{align*}

\

\section{\bf Besov classes $\bs{B_{\be,1}^1(\R^d)}$}
\setcounter{equation}{0}
\label{besov}

\

In this paper we need only Besov classes $B_{\be,1}^1(\R^d)$ of functions on the Euclidean space $\R^d$. We give here a brief introduction to such spaces and we refer the reader to \cite{Pee} for detailed information about Besov classes.

Let $w$ be an infinitely differentiable function on $\R$ such
that
\bay
\label{w}
w\ge0,\quad\supp w\subset\left[\frac12,2\right],\quad\mbox{and} \quad w(s)=1-w\left(\frac s2\right)\quad\mbox{for}\quad s\in[1,2].
\ey

We define the functions $W_n$, $n\in\Z$, on $\R^d$ by 
\bay
\label{opredWn}
\big(\F W_n\big)(x)=w\left(\frac{\|x\|_2}{2^n}\right),\quad n\in\Z, \quad x=(x_1,\cdots,x_d),
\quad\|x\|_2\df\left(\sum_{j=1}^dx_j^2\right)^{1/2},
\ey
where $\F$ is the {\it Fourier transform} defined on $L^1\big(\R^d\big)$ by
$$
\big(\F f\big)(t)=\!\int\limits_{\R^d} f(x)e^{-{\rm i}(x,t)}\,dx,\!\quad 
x=(x_1,\cdots,x_d),
\quad t=(t_1,\cdots,t_d), \!\quad(x,t)\df \sum_{j=1}^dx_jt_j.
$$
Clearly,
$$
\sum_{n\in\Z}(\F W_n)(t)=1,\quad t\in\R^d\setminus\{0\}.
$$

With each tempered distribution $f\in{\mathscr S}^\prime\big(\R^d\big)$, we
associate the sequence $\{f_n\}_{n\in\Z}$,
\bay
\label{fn}
f_n\df f*W_n.
\ey
The formal series
$
\sum_{n\in\Z}f_n
$
is a Littlewood--Paley type expansion of $f$. This series does not necessarily converge to $f$.

Initially we define the (homogeneous) Besov class $\dot B^1_{\be,1}\big(\R^d\big)$ as the space of 
$f\in{\mathscr S}^\prime(\R^n)$
such that
\bay
\label{Wn}
\{2^{n}\|f_n\|_{L^\be}\}_{n\in\Z}\in\ell^1(\Z)
\ey
and put
$$
\|f\|_{B^1_{\be,1}}\df\big\|\{2^{n}\|f_n\|_{L^\be}\}_{n\in\Z}\big\|_{\ell^1(\Z)}.
$$
According to this definition, the space $\dot B^1_{\be,1}(\R^n)$ contains all polynomials
and all polynomials $f$ satisfy the equality $\|f\|_{B^s_{p,q}}=0$. Moreover, the distribution $f$ is determined by the sequence $\{f_n\}_{n\in\Z}$
uniquely up to a polynomial. It is easy to see that the series 
$\sum_{n\ge0}f_n$ converges in ${\mathscr S}^\prime(\R^d)$.
However, the series $\sum_{n<0}f_n$ can diverge in general. It can easily be proved that the series
\bay
\label{ryad}
\sum_{n<0}\frac{\partial f_n}{\partial x_j},\quad \mbox{where}\quad 1\le j\le d,
\ey
converges uniformly on $\R^d$.


Now we can define the modified (homogeneous) Besov class $B^1_{\be,1}\big(\R^d\big)$. We say that a tempered distribution $f$
belongs to $B^1_{\be,1}(\R^d)$ if \rf{Wn} holds and
$$
\frac{\partial f}{\partial x_j}
=\sum_{n\in\Z}\frac{\partial f_n}{\partial x_j},\quad
1\le j\le d,
$$
in the space ${\mathscr S}^\prime\big(\R^d\big)$ (equipped with the weak-$*$ topology). Now the function $f$ is determined uniquely by the sequence $\{f_n\}_{n\in\Z}$ up
to a constant polynomial, and a polynomial $g$ belongs to 
$B^1_{\be,1}\big(\R^d\big)$
if and only if $g$ is constant.

Note that the functions $f_n$ defined by \rf{fn} have the following properties: $f_n\in L^\be(\R^d)$ and 
$\supp\F f\subset\{\xi\in\R^d:~\|\xi\|\le2^{n+1}\}$. Bounded continuous functions whose Fourier transforms are supported in $\{\xi\in\R^d:~\|\xi\|\le\s\}$ can be characterized by the following Paley--Wiener--Schwartz type theorem  (see \cite{R}, Theorem 7.23 and exercise 15 of Chapter 7):

\medskip

{\it Let $f$ be a continuous function
on $\R^d$ and let $M,\,\s>0$. The following statements are equivalent:

{\em(i)} $|f|\le M$ and $\supp\F f\subset\{\xi\in\R^d:\|\xi\|\le\s\}$;

{\em(ii)} $f$ is a restriction to $\R^d$ of an entire function on $\C^d$ such that 
$$
|f(z)|\le Me^{\s\|\im z\|}
$$
for all $z\in\C^d$.}

We need one more elementary remark on the Besov classes $B^1_{\be,1}(\R^d)$.

\medskip

{\bf Remark.} Suppose that $\{g_j\}_{j\ge0}$ is a sequence of functions in
$L^\be(\R^d)$ such that
$$
\quad\supp\F g_j\subset\big[-2^{j+1},2^{j+1}\big]^d\quad\mbox{and}\quad
\sum_{j\ge0}2^j\|g_j\|_{L^\be(\R^d)}<\be,
$$
then the series $\sum_{j\ge0}g_j$ converges uniformly on $\R^d$, the sum of the series belongs to $B^1_{\be,1}(\R^d)$ and
$$
\left\|\sum_{j\ge0}g_j\right\|_{B^1_{\be,1}(\R^d)}\le\const
\sum_{j\ge0}2^j\|g_j\|_{L^\be(\R^d)}.
$$

\

\section{\bf The main result}
\setcounter{equation}{0}
\label{osnre}

\

To establish the main result of the paper, we introduce the classes 
$\E^p_\s(\R^d)$, $1\le p\le\be$, $\s>0$. Put
$$
\E^p_\s(\R^d)\df\left\{f\in L^p(\R_d):~\supp\F f\subset[-\s,\s]^d\right\}.
$$

\begin{thm}
\label{Lip3p}
Let $p\in[1,\be]$.
There is no constant $K>0$ such that 
\begin{align}
\label{nerK}
\|f(A_1,B_1,C_1)&-f(A_2,B_2,C_2)\|_{\bS_p}\nonumber\\[.2cm]
&\le K\|f\|_{B_{\be,1}^1(\R^3)}
\max\big\{\|A_1-A_2\|_{\bS_p},\|B_1-B_2\|_{\bS_p},\|C_1-C_2\|_{\bS_p}\big\}
\end{align}
for all triples of not necessarily commuting finite rank self-adjoint operators $(A_1,B_1,C_1)$ and $(A_2,B_2,C_2)$ and all functions $f$ in $B_{\be,1}^1(\R^3)$.
\end{thm}

We need the following elementary lemma:

\begin{lem}
\label{proizve}
Let $\psi$ be an infinitely differentiable function on $\R$ with compact support and such that
$\psi(t)=t$ for $t\in[-1,1]$. Suppose that $\s>0$. There exists a positive number
$\vk$ such that
$$
\|\f\otimes\psi\|_{B_{\be,1}^1(\R^3)}\le\vk\|\f\|_{L^\be(\R^2)}
$$
for an arbitrary function $\f$ in $\E^\be_\s(\R^2)$, where the function 
$\f\otimes\psi$ on $\R^3$ is defined by
$$
(\f\otimes\psi)(x,y,z)\df\f(x,y)\psi(z),\quad(x,y,z)\in\R^3.
$$
\end{lem}

Let us first deduce Theorem \ref{Lip3p} from Lemma \ref{proizve} and then prove Lemma \ref{proizve}.

\medskip

{\bf Proof of Theorem \ref{Lip3p}.} Let $\psi$ be a function on $\R$ that satisfies the hypotheses of Lemma \ref{proizve}.
Let $\f$ be a function in $\E^\be_1(\R^2)$. We define the function $f$ on $\R^3$ by
$$
f(x,y,z)\df\f(x,y)\psi(z),\quad(x,y,z)\in\R^3.
$$
Suppose that $A$, $B$ and $C$ are finite rank self-adjoint operators. We consider the triples $(A,B,C)$ and $(A,B,\0)$, where $\0$ is the zero operator. It is easy to see that if $\|C\|\le1$, then $\psi(C)=C$ and
\begin{align}
\label{dvukhet}
f(A,B,C)-f(A,B,\0)&=\f(A,B)\big(\psi(C)-\psi(\0)\big)\nonumber\\[.2cm]
&=\f(A,B)\big(C-\0\big)=\f(A,B)C.
\end{align}

Let us construct the operators $A$ and $B$ and the function $\f$. The construction is similar to the construction given in the proof of Theorem 8.1 of \cite{ANP}.

Let $\{g_j\}_{1\le j\le N}$ and $\{h_j\}_{1\le j\le N}$
be orthonormal systems in Hilbert space. 
Consider the rank one projections $P_j$ and $Q_j$ defined by
$$
P_j v=(v,g_j)g_j\quad\mbox{and}\quad Q_jv=(v,h_j)h_j,\quad 1\le j\le N.
$$
We define the function $\eta$ on $\R$ by
$$
\eta(x)=\frac{2(1-\cos x)}{x^2},\quad x\in\R,\quad x\ne0,
$$
and extend it to $\R$ by continuity.
It is well known and it is easy to verify that $\eta\in\E_1^\be(\R)$. Clearly, $\eta(0)=1$ and $\eta(2k\pi)=0$,
$k\in\Z\setminus\{0\}$. Put
$$
\eta_j(x)\df\eta(x-2\pi j),\quad j\in\Z.
$$
Suppose that $\{\t_{jk}\}_{1\le j,k\le N}$ is a family of complex numbers. Define the function $\f$ by
\bay
\label{opredpsi}
\f(x,y)=\sum_{j,k}\t_{jk}\eta_j(x)\eta_k(y).
\ey
Then  $\f\in\E^\be_1(\R^2)$ and
\bay
\label{otsecjk}
\|\f\|_{L^\be(\R^2)}\le\const\max_{j,k}|\t_{jk}|,
\ey
see \cite{ANP}, \S\;8. We define now the finite rank self-adjoint operators $A$ and $B$ by
$$
A=\sum_{j=1}^n2\pi jP_j\quad\mbox{and}\quad
B=\sum_{k=1}^N2\pi kQ_k.
$$
It follows from \rf{opredpsi} that
$$
\f(A,B)=\sum_{j=1}^N\sum_{k=1}^N\f(2\pi j,2\pi k)P_jQ_k
=\sum_{j=1}^N\sum_{k=1}^N\t_{jk}P_jQ_k
=\sum_{j=1}^N\sum_{k=1}^N\t_{jk}(h_k,g_j)(\cdot\,,h_k)g_j.
$$
In other words,
$$
\f(A,B)u=\sum_{j=1}^N\sum_{k=1}^N\t_{jk}(h_k,g_j)(u,h_k)g_j
$$
for every vector $u$.

Clearly, for every unitary matrix
$\{u_{jk}\}_{1\le j,k\le N}$, there exist orthonormal systems $\{g_j\}_{1\le j\le N}$ and $\{h_j\}_{1\le j\le N}$
such that  $(h_k,g_j)=u_{jk}$.
Put 
$$
u_{jk}\df\frac1{\sqrt N}\exp\left(\frac{2\pi{\rm i}jk}N\right),\quad 1\le j,k\le N.
$$
Obviously, $\{u_{jk}\}_{1\le j,k\le N}$ is a unitary matrix.
Hence, 
we may find vectors $\{g_j\}_{j=1}^N$ and $\{h_j\}_{j=1}^N$ such that $(h_k,g_j)=u_{jk}$.
Put $\tau_{jk}=\sqrt N \,\,\ov u_{jk}$. By \rf{otsecjk},
\bay
\label{ravogr}
\|f\|_{L^\be(\R^2)}\le\const
\ey
and
$$
\f(A,B)=\sum_{j=1}^N\sum_{k=1}^N|u_{jk}|(\cdot\,,h_k)g_j
=\frac1{\sqrt N}\left(\cdot\:,\sum_{k=1}^Nh_k\right)\sum_{j=1}^Ng_j.
$$
We can define now the rank one self-adjoint operator $C$ by
\bay
\label{opopC}
C=\frac1{N}\left(\cdot\:,\sum_{j=1}^Nh_j\right)\sum_{j=1}^Nh_j.
\ey
Clearly, $\|C\|=1$ and by \rf{dvukhet},
$$
\f(A,B)C=\frac1{N^{3/2}}
\left(\sum_{j=1}^Nh_j,\sum_{j=1}^Nh_j\right)
\left(\cdot\:,\sum_{j=1}^Nh_j\right)\sum_{j=1}^Ng_j
=\frac1{\sqrt N}\left(\cdot\:,\sum_{j=1}^Nh_j\right)\sum_{j=1}^Ng_j.
$$
It is easy to see that $\f(A,B)C$ is a rank one self-adjoint operator
and
$$
\|\f(A,B)C\|_{\bS_p}=N^{1/2}
$$
for every $p\in[1,\be]$.

The result follows now from \rf{dvukhet} and \rf{ravogr}. $\bl$

\medskip

{\bf Remark.} Clearly, we can multiply the operator $C$ in \rf{opopC} by $\e_n$,
$n\ge1$, where $\{\e_n\}$ is a sequence of positive numbers such that $\e_n\le1$ and $\e_n\to0$ as $n\to\be$. This allows us to say that there are sequences 
$\{A_n\}$, $\{B_n\}$,
$\big\{C_n^{(1)}\big\}$ and $\big\{C_n^{(2)}\big\}$ of finite rank self-adjoint
operators and a sequence $\{f_n\}$ of functions in $B_{\be,1}^1(\R^3)$ such that
$$
\|f_n\|_{B^\be_{\be,1}}\le\const,
$$
$$
\lim_{n\to\be}\big\|C_n^{(1)}-C_n^{(2)}\big\|_{\bS_p}=0,
$$
but
$$
\lim_{n\to\be}
\big\|f_n\big(A_n,B_n,C_n^{(1)}\big)-f_n\big(A_n,B_n,C_n^{(2)}\big)\big\|_{\bS_p}
=\be.
$$

\medskip

{\bf Proof of Lemma \ref{proizve}.} It is well known that such functions $\psi$ belong to all Besov classes, see \cite{Pee}. Let $\psi_n\df\psi*W_n$, where the $W_n$ are defined in \rf{opredWn}. Since $\psi\in B^1_{\be,1}(\R)$, we have
$$
\sum_{n\in\Z}2^n\|\psi_n\|_{L^\be(\R)}<\be,
$$
and so
$$
\sum_{n\ge0}\|\psi_n\|_{L^\be(\R)}<\be.
$$
Put now 
$$
\psi^\flat\df\psi-\sum_{n\ge0}\psi_n.
$$
Clearly, $\psi^\flat\in L^\be(\R)$. 

It is easy to see that
$$
\supp(\f\otimes\psi^\flat)\subset[-1,1]^3\qquad\mbox{and}\qquad
\supp(\f\otimes\psi_n)\subset\Big[-2^{n+1},2^{n+1}\Big]^3,\quad n\ge0.
$$

By the remark at the end of \S\;\ref{besov}, we have
\begin{align*}
\|\f\otimes\psi\|_{B^1_{\be,1}(\R^3)}&\le
\big\|\f\otimes\psi^\flat\big\|_{B^1_{\be,1}(\R^3)}+
\left\|\f\otimes\sum_{n\ge0}\psi_n\right\|_{B^1_{\be,1}(\R^3)}\\[.2cm]
&\le\const\|\f\|_{L^\be(\R^2)}\left(\big\|\psi^\flat\big\|_{L^\be(\R)}
+\sum_{n\ge0}2^n\|\psi_n\|_{L^\be(\R)}
\right)
\end{align*}
which completes the proof. $\bl$

\

\section{\bf Lipschitz type estimates in terms of the rank of the operators}
\setcounter{equation}{0}
\label{rank}

\

In this section we consider the problem to obtain a Lipschitz type estimate for functions of finite rank noncommuting self-adjoint operators in terms of their rank.

Let us first consider the case of pairs of finite rank self-adjoint operators.
Recall that it was proved in \cite{ANP} that for $p\in[1,2]$, we have the following Lipschitz type estimate:
$$
\|f(A_1,B_1)-f(A_2,B_2)\|_{\bS_p}\le\const\|f\|_{B^1_{\be,1}(\R^2)}
\max\big\{\|A_1-A_2\|_{\bS_p},\|A_1-A_2\|_{\bS_p}\big\}
$$
for arbitrary pairs $(A_1,B_1)$ and $(A_2,B_2)$ of self-adjoint operators and for arbitrary functions $f$ in $B^1_{\be,1}(\R^2)$. On the other hand, the reasoning given in the proof of Theorem 8.1 of \cite{ANP} shows that for $p\in[2,\be]$, there exist a sequence $\{f_N\}$ of functions in $\E^\be_2(\R^2)$, sequences 
$\big\{A_1^{(N)}\big\}$, $\big\{A_2^{(N)}\big\}$ and
$\big\{B^{(N)}\big\}$ of self-adjoint operators of rank at most $N$ such that
$$
\|f_N\|_{L^\be(\R^2)}\le\const,
$$
and
$$
\big\|\big(f(A_1^{(N)},B^{(N)}\big)-\big(f(A_2^{(N)},B^{(N)}\big)\big\|_{\bS_p}
\ge\const N^{1/2-1/p}\big\|A_1^{(N)}-A_2^{(N)}\big\|_{\bS_p}.
$$

The following result shows that this estimate is sharp.

\begin{thm}
\label{luneby}
Let $(A_1,B_1)$ and $(A_2,B_2)$ be pairs of self-adjoint operators of rank at most $N$ and let $p\in[2,\be]$. Then
$$
\big\|\big(f(A_1,B_1\big)-\big(f(A_2,B_2\big)\big\|_{\bS_p}
\le\const N^{1/2-1/p}\;
\|f\|_{B^1_{\be,1}}\max\big\{\|A_1-A_2\|_{\bS_p},\|B_1-B_2\|_{\bS_p}\big\}
$$
for every function $f$ in $B_{\be,1}^1(\R^2)$
\end{thm}

\Pf By Theorem 7.2 of \cite{ANP},
$$
\big\|\big(f(A_1,B_1\big)-\big(f(A_2,B_2\big)\big\|_{\bS_2}
\le\const 
\|f\|_{B^1_{\be,1}}\max\big\{\|A_1-A_2\|_{\bS_2},\|B_1-B_2\|_{\bS_2}\big\}.
$$
Obviously,
$$
\big\|\big(f(A_1,B_1\big)-\big(f(A_2,B_2\big)\big\|_{\bS_p}
\le\big\|\big(f(A_1,B_1\big)-\big(f(A_2,B_2\big)\big\|_{\bS_2}.
$$
The result follow from the following well-known inequalities for finite rank operators:
$$
\|A_1-A_2\|_{\bS_2}\le\const N^{1/2-1/p}\|A_1-A_2\|_{\bS_p}
$$
and
$$
\|B_1-B_2\|_{\bS_2}\le\const N^{1/2-1/p}\|B_1-B_2\|_{\bS_p}.\quad\bl
$$

A similar problem can be posed in the case of functions of triples of not necessarily commuting self-adjoint operators of finite rank. The reasoning given in the proof of Theorem \ref{Lip3p} shows that for $p\in[1,\be]$, there exist
a sequence $\{f_N\}$ of functions in $B_{\be,1}^1(\R^3)$, sequences 
$\big\{A^{(N)}\big\}$, $\big\{B^{(N)}\big\}$, $\big\{C_1^{(N)}\big\}$,
and $\big\{C_2^{(N)}\big\}$ of self-adjoint operators 
of rank at most $N$ such that
$$
\|f_N\|_{B_{\be,1}^1}\le\const,
$$
and
$$
\big\|\big(f(A^{(N)},B^{(N)},C_1^{(N)}\big)-
\big(f(A_2^{(N)},B_2^{(N)},C_2^{(N)}\big)\big\|_{\bS_p}
\ge\const N^{1/2}\|C_1-C_2\|_{\bS_p}.
$$

{\it I do not know whether this lower estimate is sharp}. To obtain a trivial upper estimate, we need the following elementary formula:
\begin{align}
\label{periztr}
f&(A_1,B,C)-f(A_2,B,C)\nonumber\\[.2cm]
&=\sum
\frac{f(\l_1,\mu,\nu)-f(\l_2,\mu,\nu)}{\l_1-\l_2}
E_{A_1}(\{\l_1\})(A_1-A_2)E_{A_2}(\{\l_2\})E_B(\{\mu\})E_C(\{\nu\})
\end{align}
for an arbitrary function $f$ on $\R^3$ and arbitrary finite rank self-adjoint operators
$A_1$, $A_2$, $B$ and $C$, where $E_{A_1}$, $E_{A_2}$, $E_B$ and $E_C$ are the spectral projections of $A_1$, $A_2$, $B$ and $C$
and the sum is taken over $\l_1$, $\l_2$, $\mu$, $\nu$ in $\R$ such that $\l_1\ne\l_2$. Formula \rf{periztr}
can be proved elementarily.

Similar formulae hold for
the differences $f(A,B_1,C)-f(A,B_2,C)$ and $f(A,B,C_1)-f(A,B,C_2)$.

Such formulae imply the following trivial upper estimate {\it for arbitrary Lipschitz functions} on $\R^3$:

\begin{thm}
\label{trivero}
Let $f$ be a Lipschitz function on $\R^3$. Suppose that $A_1$, $B_1$, $C_1$,
$A_2$, $B_2$ and $C_2$ are self-adjoint operators of rank at most $N$. Then
for $p\in[1,\be]$, the following estimate holds:
$$
\|f(A_1,B_1,C_1)-f(A_2,B_2,C_2)\|_{\bS_p}\le N^4\|f\|_{\Li}
\big(\|A_1-A_2\|_{\bS_p}+\|B_1-B_2\|_{\bS_p}+\|C_1-C_2\|_{\bS_p}\big).
$$
\end{thm}

\Pf It follows immediately from formula \rf{periztr} that 
$$
\|f(A_1,B_1,C_1)-f(A_2,B_1,C_1)\|_{\bS_p}
\le N^4\|f\|_{\Li}\|A_1-A_2\|_{\bS_p}.
$$
In the same way one can establish the inequalities:
$$
\|f(A_2,B_1,C_1)-f(A_2,B_2,C_1)\|_{\bS_p}
\le N^4\|f\|_{\Li}\|B_1-B_2\|_{\bS_p}
$$
and
$$
\|f(A_2,B_2,C_1)-f(A_2,B_2,C_2)\|_{\bS_p}
\le N^4\|f\|_{\Li}\|C_1-C_2\|_{\bS_p}
$$
which proves the result. $\bl$

\

\

\noindent
Department of Mathematics\\
Michigan State University\\
East Lansing Michigan 48824\\
and\\ 
RUDN University, 6 Miklukho-Maklay St., \\
Moscow, 117198, Russia


\begin{thebibliography}{99}
\label{bibl}


%

%
%
%

\bibitem[AP]{AP}{\sc A.B. Aleksandrov} and {\sc V.V. Peller}, {\it Operator Lipschitz functions}, Uspekhi Matem. Nauk. {\bf71:4} (2016), 3--106 (Russian).

English transl.: Russian Math. Surveys, {\bf71:4} (2016), 605--702.

\bibitem[ANP]{ANP}{\sc A.B. Aleksandrov, F.L. Nazarov} and
{\sc V.V. Peller}, {\em Functions of noncommuting self-adjoint operators under perturbation and estimates of triple operator integrals}, Adv. Math. {\bf295} (2016), 1�-52.



\bibitem[APPS]{APPS} {\sc A.B. Aleksandrov, V.V. Peller, D. Potapov}, and
{\sc F. Sukochev}, {\em Functions of normal operators under perturbations},
Advances in Math. {\bf226} (2011), 5216-�5251.


\bibitem[BS1]{BS1} {\sc M.S. Birman} and {\sc M.Z. Solomyak},
{\em Double Stieltjes operator integrals},
Problems of Math. Phys., Leningrad. Univ. {\bf1} (1966), 33--67 (Russian).

English transl.: Topics Math. Physics {\bf1} (1967), 25--54, Consultants Bureau Plenum
Publishing Corporation, New York.

\bibitem[BS2]{BS2} {\sc M.S. Birman} and {\sc M.Z. Solomyak},
 {\em Double Stieltjes operator integrals. II},
 Problems of Math. Phys., Leningrad. Univ. {\bf2} (1967), 26--60 (Russian).

English transl.: Topics Math. Physics {\bf2} (1968), 19--46, Consultants Bureau Plenum
Publishing Corporation, New York.

%

\bibitem[BS3]{BS3} {\sc M.S. Birman} and {\sc M.Z. Solomyak},
{\em Double Stieltjes operator integrals. III},
Problems of Math. Phys., Leningrad. Univ. {\bf6} (1973), 27--53 (Russian).




\bibitem[DK]{DK} {\sc Yu.L. Daletskii} and {\sc S.G. Krein}, {\em Integration and differentiation of
functions of Hermitian operators and application to the theory of perturbations} (Russian), Trudy Sem.
Functsion. Anal., Voronezh. Gos. Univ. {\bf1} (1956), 81--105.


\bibitem[F1]{F1}  {\sc Yu.B. Farforovskaya}, {\em  The connection of the Kantorovich-Rubinshtein metric for spectral resolutions of selfadjoint operators with functions of operators},
Vestnik Leningrad. Univ.  {\bf19}  (1968), 94--97. (Russian).

\bibitem[F2]{F2}  {\sc Yu.B. Farforovskaya}, {\em An example of a Lipschitzian function of selfadjoint
operators that yields a nonnuclear increase under a nuclear perturbation}.  Zap. Nauchn. Sem.
Leningrad. Otdel. Mat. Inst. Steklov. (LOMI)  {\bf30}  (1972), 146--153 (Russian).



%

%

%

\bibitem[JTT]{JTT} {\sc K. Juschenko, I.G. Todorov} and {\sc L. Turowska}, {\em Multidimensional operator multipliers}, Trans. Amer. Math. Soc. {\bf361}
(2009), 4683-�4720.

%

\bibitem[Ka]{Ka}  {\sc T. Kato}, {\em Continuity of the map $S\mapsto \mid S\mid $ for linear operators},
Proc. Japan Acad.  {\bf49}  (1973), 157--160.


\bibitem[KPSS]{KPSS} {\sc E. Kissin, D. Potapov, V. S. Shulman} and {\sc F. Sukochev}, {\it Operator smoothness in Schatten norms for functions of several variables: Lipschitz conditions, differentiability and unbounded derivations},  Proc. Lond. Math. Soc. (3) {\bf105} (2012), 661--702.





\bibitem[Mc]{Mc}
{\sc A. McIntosh}, {\em Counterexample to a question on commutators},
Proc. Amer. Math. Soc. {\bf29} (1971) 337--340.

\bibitem[NP]{NP} {\sc F.L. Nazarov} and {V.V. Peller} {\em Functions of $n$-tuples of commuting self-adjoint operators}, J. Funct. Anal. {\bf266} (2014), 5398�-5428.


\bibitem[Pee]{Pee} {\sc J. Peetre},
{\em New thoughts on Besov spaces}, Duke Univ. Press., Durham, NC, 1976.


\bibitem[Pe1]{Pe1} {\sc V.V.Peller}, {\em Hankel operators of class ${\bf S}_{p}$
and their applications (rational approximation, Gaussian processes,
the problem of majorizing operators)}, Mat. Sbornik,
{\bf 113} (1980), 538-581.

English Transl. in Math. USSR Sbornik, {\bf 41}
(1982), 443-479.


\bibitem[Pe2]{Pe2} {\sc V.V. Peller},
{\em Hankel operators in the theory of perturbations of unitary and self-adjoint operators},
Funktsional. Anal. i Prilozhen. {\bf19:2}  (1985),
37--51 (Russian).

English transl.: Funct. Anal. Appl. {\bf19} (1985) , 111--123.



\bibitem[Pe3]{Pe3} {\sc V.V. Peller},
{\em Hankel operators in the perturbation theory of of unbounded self-adjoint operators}.
Analysis and partial differential equations,  529--544,
Lecture Notes in Pure and Appl. Math., {\bf122}, Dekker, New York, 1990.

\bibitem[Pe4]{Pe4} {\sc V.V. Peller}, {\em Hankel operators and their applications,}
Springer-Verlag, New York, 2003.

\bibitem[Pe5]{Pe5} {\sc V.V. Peller}, {\em Multiple operator integrals and higher operator
derivatives}, J. Funct. Anal.  {\bf233}  (2006),  515--544.

\bibitem[Pe6]{Pe6} {\sc V.V. Peller}, {\em Multiple operator integrals in perturbation theory}, Bull. Math. Sci. {\bf6} (2016), 15--88.

\bibitem[Pe7]{Pe7} {\sc V.V. Peller}, {Functions of triples of noncommuting self-adjoint operators and their perturbations},  arXiv:1606.0896.


\bibitem[Pi]{Pi} {\sc G. Pisier}, {\em Similarity problems and completely bounded maps},
Second, expanded edition. Includes the solution to ``The Halmos problem''. Lecture Notes in Mathematics,
1618. Springer-Verlag, Berlin, 2001.

\bibitem[PS]{PS} {\sc D. Potapov} and {\sc F. Sukochev}, {\em Operator-Lipschitz functions in Schatten--von Neumann classes}, Acta Math. {\bf207} (2011), 375--389.

\bibitem[R]{R}{\sc W. Rudin}, {\em Functional analysis}, M$^{\rm c}$Graw Hill, 1991.


%

%









%
%

\end{thebibliography}
\end{document}